# ESTIMATING MARGINAL SURVIVAL FUNCTION BY ADJUSTING FOR DEPENDENT CENSORING USING MANY COVARIATES[1]

By Donglin Zeng

*University of North Carolina at Chapel Hill*

One goal in survival analysis of right-censored data is to estimate the marginal survival function in the presence of dependent censoring. When many auxiliary covariates are sufficient to explain the dependent censoring, estimation based on either a semiparametric model or a nonparametric model of the conditional survival function can be problematic due to the high dimensionality of the auxiliary information. In this paper, we use two working models to condense these high-dimensional covariates in dimension reduction; then an estimate of the marginal survival function can be derived nonparametrically in a low-dimensional space. We show that such an estimator has the following double robust property: when either working model is correct, the estimator is consistent and asymptotically Gaussian; when both working models are correct, the asymptotic variance attains the efficiency bound.

**1. Introduction.** Right-censored data with dependent censoring are common in many epidemiological studies. Such data consist of $n$ i.i.d. copies of the observation $(Y = T \wedge C, R = I(T \leq C), L)$, where $T$ is the failure time of interest, $C$ is the right censoring time, and $L$ includes the covariate information. Usually, the covariates $L$ contain not only subject demographic information and disease history, but also much other auxiliary information which researchers are not primarily interested in but which is informative in predicting subjects' failure time or explaining why subjects drop out, or both. For example, in a typical medical study, $L$ may contain the patient's willingness to participate in the study, the patient's accessibility to hospitals, the social support from the patient's family members, or the patient's

Received September 2001; revised September 2003.
[1]Supported in part by NIDA Grant 2P50 DA 10075 and NSF Grant DMS-98-02885.
*AMS 2000 subject classifications.* Primary 62G07; secondary 62F12.
*Key words and phrases.* Dependent censoring, double robustness, semiparametric efficiency.

genetic information, and so on. When much auxiliary information has been collected, in practice, it is safe to assume that $L$ is sufficient to explain the dependence between $T$ and $C$. Equivalently, $T$ and $C$ are independent when conditional on $L$.

The purpose of this article is to estimate the marginal survival function of $T$ using right-censored data. A standard estimate is the Kaplan–Meier estimate. However, it is well known that, when $T$ and $C$ are dependent, this estimator is inconsistent. Another intuitive approach to estimate the survival function of $T$ is to estimate the conditional distribution of $T$ given $L$ using a semiparametric model [e.g., the Cox proportional hazard model; see Cox (1972)], the proportional odds model [Bennett (1983), etc.], or via nonparametric estimation approaches such as using a local likelihood function [Tibshirani and Hastie (1987)]. Then the estimate of the marginal survival function of $T$ is simply the empirical average of the conditional distribution of $T$ given $L$ over all the observed covariates. However, the above approaches can be problematic when the auxiliary covariates, $L$, consist of many variables. This is because when $L$ has at least three dimensions, nonparametric estimation of the distribution of $T$ given $L$ is infeasible in a moderate-sized sample due to the curse of dimensionality; and in any semiparametric model, the parametric function of $L$ in the model of $T$ given $L$ is likely to be misspecified. Consequently, these intuitive approaches bias the estimation of the survival function of $T$.

To reduce the limitation in the above intuitive approaches, in this article we propose two working models for both the lifetime $T$ and the censoring time $C$ given all the covariates $L$. Then two-dimensional condensed information of $L$ is extracted from the working models and used as the new covariates in place of $L$. The estimator of the survival function is obtained by maximizing a pseudo-likelihood function nonparametrically in the space with the reduced dimension. It is shown that if either working model is correct, the estimator of the marginal survival function is consistent and asymptotically Gaussian; if both working models are correct, the asymptotic variance of the estimator attains the generalized Cramér–Rao bound of the full model space [cf. Bickel, Klaassen, Ritov and Wellner (1993)]. The first property is named "double robustness" by Robins, Rotnitzky and van der Laan (2000), since the estimator remains consistent if one working model is misspecified but the other one is correct.

The method of using the condensed information of the high-dimensional covariates in the estimation dates back to the propensity score approach by Rubin (1976) in a simple regression, where the propensity score was defined as the predicted missing probability given all the covariates. Little (1986) further combined the propensity score and the mean score, the latter of which was defined as the predicted mean response given all the covariates, to estimate the population mean in a survey study. Such methods have been recently developed and generalized to study dependent censoring in semiparametric regression and survival analysis by Robins and others [Rotnitzky and Robins (1995), Robins, Rotnitzky and van der Laan (2000) and Scharfstein and Robins (2002)]. Although all the above mentioned approaches including ours pursue the summary information of the covariates, sometimes referred to as the propensity score or risk score, using the working models for $T$ and $C$ given $L$, the estimation approach we take is much different from theirs. Robins, Rotnitzky and van der Laan's approach is to begin with an inverse-weighted estimating equation, where only the complete observations are used in the estimating equation and each complete observation is weighted with the inverse of the probability of not being censored; a final estimating equation for estimating the marginal survival distribution is to subtract from the inverse-weighted estimating equation the projection on the score tangent space. However, the method we propose in this article is a purely likelihood-based approach: we first obtain the condensed covariates by optimizing the pseudolikelihood functions based on the working models; we then optimize another pseudolikelihood function to derive the estimate of $T$'s survival function. Therefore, the likelihood-based approach we take only involves simple optimization steps and the estimate turns out to have a simple expression; by contrast, the approach in Robins, Rotnitzky and van der Laan (2000) requires a practical user to have knowledge of the projection on the score space.

This article is organized as follows: In Section 2 we give the details of estimating the marginal survival function; the asymptotic properties of our estimator are then given in Section 3, where we also provide an algorithm to estimate the asymptotic variance; the numerical results from a simulation study are given in Section 4; finally, the article concludes with some discussion. Most of the proofs in this article are deferred to the Appendix.

**2. Estimation.** Under the assumption that $T$ and $C$ are independent given $L$, the observed likelihood function for $n$ observations can be written as

$$\prod_{i=1}^{n} [h_{T|L}(Y_i|L_i)^{R_i} e^{-H_{T|L}(Y_i|L_i)} h_{C|L}(Y_i|L_i)^{1-R_i} e^{-H_{C|L}(Y_i|L_i)} f_L(L_i)],$$

where $h_{T|L}(\cdot|L)$ and $h_{C|L}(\cdot|L)$ are the hazard rate functions for $T$ and $C$ given $L$, respectively; $H_{T|L}(\cdot|L)$ and $H_{C|L}(\cdot|L)$ are their respective cumulative hazard functions. Our estimation procedure consists of the following steps.

*Step* 1. We propose two working models for both the lifetime $T$ and the censoring time $C$ given $L$. Our working models for $T$ given $L$ and $C$ given

$L$ are Cox's proportional hazard models; that is, we tentatively assume that

$$h_{T|L}(y|l) = \lambda_T(y)e^{\beta' l}, \qquad h_{C|L}(y|l) = \lambda_C(y)e^{\gamma' l}$$

for some unknown functions $\lambda_T(\cdot), \lambda_C(\cdot)$ and some parameters $(\beta, \gamma)$.

*Step* 2. We derive the estimator of $(\beta, \gamma)$ simply by performing Cox's regressions, or equivalently, we maximize the following pseudolog partial likelihood functions:

$$\tilde{L}_1^{(n)}(\beta) = \frac{1}{n}\sum_{i=1}^n R_i\left[\beta' L_i - \log\left(\sum_{Y_j \geq Y_i} e^{\beta' L_j}\right)\right],$$

$$\tilde{L}_2^{(n)}(\gamma) = \frac{1}{n}\sum_{i=1}^n (1-R_i)\left[\gamma' L_i - \log\left(\sum_{Y_j \geq Y_i} e^{\gamma' L_j}\right)\right],$$

to estimate $\beta$ and $\gamma$, respectively. We denote the estimators as $(\hat{\beta}_n, \hat{\gamma}_n)$. It will be shown in the next section that there exist two constants $\beta^*$ and $\gamma^*$ such that $\hat{\beta}_n$ and $\hat{\gamma}_n$ converge to $\beta^*$ and $\gamma^*$ in probability, respectively.

*Step* 3. Acting as if the two limit constants $\beta^*$ and $\gamma^*$ were known, we obtain the estimator of the hazard rate function of $T$ given $(\beta^{*'}L, \gamma^{*'}L)$ as follows. Denote $Z^* = (\beta^{*'}L, \gamma^{*'}L)$. When either of the working models is right, it will be shown that $T$ and $C$ are independent given $Z^*$ in Lemma 3.1. In other words, the two-dimensional covariate $Z^*$ is sufficient to explain the dependence between $T$ and $C$. Therefore, we replace the covariates $L$ by $Z^*$ in the observations and obtain a reduced dataset $(Y_i, R_i, Z_i^* = (\beta^{*'}L_i, \gamma^{*'}L_i)), i = 1, \ldots, n$. Clearly, the likelihood function for this reduced data can be verified to be

$$\prod_{i=1}^n [h_{T|Z^*}(Y_i|Z_i^*)^{R_i} e^{-H_{T|Z^*}(Y_i|Z_i^*)} h_{C|Z^*}(Y_i|Z_i^*)^{1-R_i} e^{-H_{C|Z^*}(Y_i|Z_i^*)} f_{Z^*}(Z_i^*)],$$

where $h_{T|Z^*}(\cdot|Z^*), h_{C|Z^*}(\cdot|Z^*)$ are the hazard rate functions of $T$ and $C$ given $Z^*$, respectively, and $H_{T|Z^*}(\cdot|Z^*), H_{C|Z^*}(\cdot|Z^*)$ are their corresponding cumulative hazard functions. So we can estimate $h_{T|Z^*}(y|z)$ by maximizing a local version of the observed log-likelihood function

$$\sum_{i=1}^n K\left(\frac{Z_i^* - z}{a_n}\right)[R_i \log h_{T|Z^*}(Y_i|z) - H_{T|Z^*}(Y_i|z)],$$

where $K(\cdot, \cdot)$ is a symmetric two-dimensional kernel function and $a_n$ is a bandwidth to be chosen later. Easy calculation shows that the maximizer for $h_{T|Z^*}(y|z)$ is an empirical function with a point mass at each observed $Y_j$ and the mass is equal to $R_j K(\frac{Z_j^* - z}{a_n})/(\sum_{Y_m \geq Y_j} K(\frac{Z_m^* - z}{a_n}))$.

*Step* 4. Therefore, the estimator for the cumulative hazard function is given by

$$\hat{H}_{T|Z^*}(y|z) = \sum_{Y_j \leq y} \frac{R_j K((Z_j^* - z)/a_n)}{\sum_{Y_m \geq Y_j} K((Z_m^* - z)/a_n)}.$$

The estimator for the conditional survival function of $T$ given $Z^*$ is then $\hat{S}_{T|Z^*}(t|z) = \prod_{s \leq t}(1 - \hat{H}_{T|Z^*}(\{s\}|z))$. Finally, the estimator for the marginal survival function of $T$ is simply the empirical average of $\hat{S}_{T|Z^*}(t|z)$ over all the $Z_i^*, i = 1, \ldots, n$. That is, it is equal to

$$\frac{1}{n} \sum_{i=1}^n \prod_{j=1}^n \left(1 - \frac{K((Z_i^* - Z_j^*)/a_n) I_{Y_j \leq t} R_j}{\sum_{m=1}^n K((Z_i^* - Z_m^*)/a_n) I_{Y_j \leq Y_m}}\right).$$

*Step* 5. Since the two constants $\beta^*$ and $\gamma^*$ are unknown but can be consistently estimated by $\hat{\beta}_n$ and $\hat{\gamma}_n$, we replace $(\beta^*, \gamma^*)$ with $(\hat{\beta}_n, \hat{\gamma}_n)$ in the last estimator obtained in Step 4. Thus, we obtain an estimator for the survival function of $T$ as

$$\hat{S}_n(t) = \frac{1}{n} \sum_{i=1}^n \prod_{j=1}^n \left(1 - \frac{K((\hat{Z}_i - \hat{Z}_j)/a_n) I_{Y_j \leq t} R_j}{\sum_{m=1}^n K((\hat{Z}_i - \hat{Z}_m)/a_n) I_{Y_j \leq Y_m}}\right).$$

**3. Main results.** Before we present the main results of this article, we assume the following conditions hold.

ASSUMPTION 3.1. *$T$ and $C$ are independent conditional on $L$.*

ASSUMPTION 3.2. *Let $\tau$ be the ending time of the study. For any $l$ in the support of $L$, the conditional density of $(T, C)$ given $L = l$ is continuously twice-differentiable in $[0, \infty) \times [0, \tau)$ and its second derivatives are uniformly bounded. Moreover, $L$ has bounded second derivative in its support.*

ASSUMPTION 3.3. *There exists an unknown constant $\theta$ such that for any $l$ in the support of $L$,*

$$\inf_l P(T \geq \tau | L = l) > \theta > 0,$$

$$\inf_l P(C \geq \tau | L = l) = \inf_l P(C = \tau | L = l) > \theta > 0 \quad \text{a.s.}$$

ASSUMPTION 3.4. *The kernel function $K(x_1, x_2)$ is continuously twice differentiable with bounded second derivatives. Moreover, it satisfies*

$$K(-x_1, -x_2) = K(x_1, x_2),$$

$$|\nabla_{x_j} K(x_1, x_2)| \leq \frac{O(1)}{1 + x_1^2 + x_2^2}, \quad j = 1, 2.$$

ASSUMPTION 3.5. $\frac{(\log a_n)^2}{na_n^2} \to 0, na_n^2 \to \infty, na_n^4 \to 0$.

REMARK 3.1. Assumption 3.3 implies that all the subjects surviving until $\tau$ will be right-censored at $\tau$, due to the end of the study. In Assumption 3.4, an example of kernel functions satisfying the conditions is $k(x_1, x_2) = \exp\{-(x_1^2 + x_2^2)\}$ or any symmetric smooth function with bounded support. The conditions in Assumption 3.5 stipulate the choice of the bandwidth and control the asymptotic bias of $\hat{S}_n(t)$ resulting from the kernel estimation. First, based on Dabrowska (1987), $(\log a_n)^2/(na_n^2) \to 0$ ensures the unform convergence of $\hat{H}_{T|Z^*}(t|z)$, a type of kernel estimator for the cumulative hazard function. Second, it is known that for a kernel smoothing estimator with bandwidth $a_n$ in the two-dimensional real space, the convergence rate is of the order $\sqrt{na_n^2}$ and the bias is of the order $a_n^2$. Such bias carries into the estimator $\hat{S}_n(t)$. Thus, $na_n^4 \to 0$ in Assumption 3.5 ensures that the asymptotic bias of $\sqrt{n}(\hat{S}_n(t) - S_0(t))$ resulting from the kernel estimation will be zero. Clearly, one choice of the bandwidth $a_n$ in Assumption 3.5 can be $O(1)n^{-\alpha}$ where $\alpha \in (\frac{1}{4}, \frac{1}{2})$ and we will use $a_n = O(n^{-1/3})$ in the subsequent simulation study.

3.1. *Asymptotic properties of $\hat{\beta}_n$ and $\hat{\gamma}_n$.*

THEOREM 3.1. *Under Assumptions 3.1–3.5, there exist $\beta^*$ and $\gamma^*$ such that*

$$\sqrt{n}(\hat{\beta}_n - \beta^*) = \frac{1}{\sqrt{n}} \sum_{i=1}^{n} S_\beta(\beta^*, Y_i, R_i, L_i) + o_p(1),$$

$$\sqrt{n}(\hat{\gamma}_n - \gamma^*) = \frac{1}{\sqrt{n}} \sum_{i=1}^{n} S_\gamma(\gamma^*, Y_i, R_i, L_i) + o_p(1)$$

*for some influence functions $S_\beta$ and $S_\gamma$. Thus, both $\sqrt{n}(\hat{\beta}_n - \beta^*)$ and $\sqrt{n}(\hat{\gamma}_n - \gamma^*)$ converge weakly to some multinormal distributions.*

Theorem 3.1 shows that $(\hat{\beta}_n, \hat{\gamma}_n)$ converges to some constants even though using Cox's proportional hazard models as working models may be wrong. Obviously, if the model of $T$ given $L$ is a Cox's proportional hazard model, then $\beta^*$ is the correct coefficient of $L$ specified in this model; if the model of $C$ given $L$ is a Cox's proportional hazard model, then $\gamma^*$ is the correct coefficient of $L$ specified in this model. Furthermore, we show that, when either working model is correct, the condensed variables $(\beta^{*\prime}L, \gamma^{*\prime}L)$ are sufficient to explain the dependence between the lifetime and the censoring time.

LEMMA 3.1. *Suppose either of the working models is right, that is, either the model for $T$ given $L$ is a Cox's proportional hazard model or the model for $C$ given $L$ is a Cox's proportional hazard model. Let $Z^* = (\beta^{*'}L, \gamma^{*'}L)$. Then $T$ and $C$ are independent given $Z^*$, and moreover, the cumulative hazard function of $T$ given $Z^*$ is equal to $\int_0^t \frac{d_u P(T \wedge C \leq u, R=1|Z^*=z)}{P(T \wedge C \geq u|Z^*=z)}$.*

PROOF. We only show that the results are true if the working model for $C$ given $L$ is a Cox's proportional hazard model. For any $t_1, t_2 > 0$,

$$\begin{aligned} P(T < t_1, C < t_2 | Z^*) &= E_{L|Z^*}[P(T < t_1|L)P(C < t_2|L)] \\ &= E_{L|Z^*}[P(T < t_1|L)P(C < t_2|Z^*)] \\ &= P(C < t_2|Z^*)P(T < t_1|Z^*). \end{aligned}$$

Therefore, $T$ and $C$ are independent conditional on $Z^*$. Hence,

$$\begin{aligned} &\int_0^t \frac{d_u P(T \wedge C \leq u, R=1|Z^*=z)}{P(T \wedge C \geq u|Z^*=z)} \\ &= -\int_0^t \frac{\frac{d}{du}\int_u^\infty [P(c \geq T|Z^*=z) - P(u \geq T|Z^*=z)]\, dF_C(c)}{P(C \geq u, T \geq u|Z^*=z)} du \\ &= \int_0^t \frac{d_u P(T \leq u|Z^*=z)}{P(T \geq u|Z^*=z)} = H_{T|Z^*}(t|z). \end{aligned} \qquad \square$$

3.2. *Asymptotic properties of the estimator $\hat{S}_n(t)$.* The main result is the asymptotic property for $\hat{S}_n(t)$ given below.

THEOREM 3.2. *Under Assumptions 3.1–3.5, if either of the two working models is correct, that is, either the model for $T$ given $L$ is a Cox's proportional hazard model or the model for $C$ given $L$ is a Cox's proportional hazard model,*

$$\sqrt{n}(\hat{S}_n(t) - S(t)) \Longrightarrow G(t) \qquad in\ l^\infty([0,\tau]),$$

*where $G(\cdot)$ is a Gaussian process.*

REMARK 3.2. Indeed, the covariance of $G(\cdot)$ has an explicit form. From the proof of Theorem 3.2, $\hat{S}_n(t)$ is an asymptotic linear estimator of $S(t)$ and its influence function, denoted as $\mathcal{A}(t; Y, R, L)$, is equal to

$$\begin{aligned} &e^{-H_{T|Z^*}(t|Z^*)} - S(t) - RI_{Y \leq t} e^{H_{T|Z^*}(Y|Z^*) + H_{C|Z^*}(Y|Z^*) - H_{T|Z^*}(t|Z^*)} \\ (3.1) \quad &+ \int_0^{t \wedge Y} e^{H_{T|Z^*}(u|Z^*) + H_{C|Z^*}(u|Z^*) - H_{T|Z^*}(t|Z^*)}\, d_u H_{T|Z^*}(u|Z^*) \\ &+ \mathcal{B}_1(t; Y, R, L) + \mathcal{B}_2(t; Y, R, L), \end{aligned}$$

where $\mathcal{B}_1(t; Y, R, L)$ is

$$-E\left[e^{-H_{T|Z^*}(t|Z^*)}\nabla_\gamma|_{\gamma=\gamma^*}\int_0^t \frac{d_u P(Y \wedge C \leq u, R=1|\gamma'L, \beta^{*'}L)}{P(T \wedge C \geq u|\gamma'L, \beta^{*'}L)}\right]$$
$$\times S_\gamma(\gamma^*, Y, R, L)$$

and $\mathcal{B}_2(t; Y, R, L)$ is

$$-E\left[e^{-H_{T|Z^*}(t|Z^*)}\nabla_\beta|_{\beta=\beta^*}\int_0^t \frac{d_u P(T \wedge C \leq u, R=1|\gamma^{*'}L, \beta'L)}{P(T \wedge C \geq u|\gamma^{*'}L, \beta'L)}\right]$$
$$\times S_\beta(\beta^*, Y, R, L).$$

Therefore, the covariance function, denoted by $r(s,t)$, for the limit Gaussian process is equal to $\mathrm{Cov}(\mathcal{A}(s; Y, R, L), \mathcal{A}(t; Y, R, L))$. Interestingly, the covariance of the limiting process $G(\cdot)$ does not depend on the choice of the kernel function or the choice of the bandwidth in deriving the estimator $\hat{S}_n(t)$.

In the expression of (3.1), the two terms $\mathcal{B}_1(t; Y, R, L)$ and $\mathcal{B}_2(t; Y, R, L)$ contribute to the variation in estimating $\hat{S}_n(t)$ due to the estimation of $\beta^*$ and $\gamma^*$. Moreover, if the working model of $T$ given $L$ is correct, by repeating the arguments in proving Lemma 3.1, we easily obtain that for any $\gamma$,

$$\int_0^t \frac{d_u P(T \wedge C \leq u, R=1|\gamma'L, \beta^{*'}L)}{P(T \wedge C \geq u|\gamma'L, \beta^{*'}L)} = H_{T|\beta^{*'}L, \gamma'L}(t|\beta^{*'}L, \gamma'L).$$

Therefore,

$$-E\left[e^{-H_{T|Z^*}(t|Z^*)}\nabla_\gamma|_{\gamma=\gamma^*}\int_0^t \frac{d_u P(T \wedge C \leq u, R=1|\beta^{*'}L, \gamma'L)}{P(T \wedge C \geq u|\beta^{*'}L, \gamma'L)}\right]$$
$$= \nabla_\gamma|_{\gamma=\gamma^*} E e^{-H_{T|\beta^{*'}L, \gamma'L}(t|\beta^{*'}L, \gamma'L)} = \nabla_\gamma|_{\gamma=\gamma^*} S(t) = 0.$$

Hence, we conclude that $\mathcal{B}_1(t; Y, R, L)$ is zero. Similarly, $\mathcal{B}_2(t; Y, R, L)$ is zero if the working model for $C$ given $L$ is correct.

COROLLARY 3.1. *In the expression of* (3.1), *if the working model for $T$ given $L$ is correct, $\mathcal{B}_1(t; Y, R, L) = 0$; if the working model for $C$ given $L$ is correct, $\mathcal{B}_2(t; Y, R, L) = 0$.*

As a result, when both working models are correct, $\mathcal{B}_1(t; Y, R, L) = \mathcal{B}_2(t; Y, R, L) = 0$ and moreover, $H_{T|Z^*}(t|Z^*) = H_{T|L}(t|L)$, $H_{C|Z^*}(t|Z^*) = H_{C|L}(t|L)$. Hence, simple calculation gives that the influence function in (3.1) for $\hat{S}_n(t)$ is equal

to

$$S_{T|L}(t|L) - S(t) + \frac{R(I(T \geq t) - S(t))}{S_{C|L}(T|L)}$$
$$+ \int E\left[\frac{R(I(T \geq t) - S(t))}{S_{C|L}(T|L)}\Big| L, T \geq u, C \geq u\right] dM_C(u),$$

where $dM_C(u) = (1-R)\, dI(Y \leq u) - I(Y \geq u)\, dH_{C|L}(u|L)$ is the martingale process for the censoring time. This turns out to be the efficient influence function for $S(t)$ in the full model space, which was derived in an unpublished manuscript by Gill, van der Laan and Robins (1997). Consequently, we have obtained the following corollary.

COROLLARY 3.2. *When both working models are correct, the asymptotic variance of $\hat{S}_n(t)$ is the same as the generalized Cramér–Rao bound for $S(t)$.*

3.3. *Variance estimation for estimating a Fréchet differentiable functional of $S(t)$.* In survival analysis, practical interest may include the estimation of some functional of $S(t)$, such as the survival probability at a fixed time $t_0$, the observed mean lifetime $E[T|T \leq \tau]$, and median lifetime, and so on. Denote such a functional of $S(t)$ as $\Psi(S(t))$. Then we can estimate it with $\Psi(\hat{S}_n(t))$. Furthermore, if $\Psi(\cdot)$ is Fréchet differentiable with its first derivative along direction $\hat{S}_n(t) - S(t)$ given by $\int_0^\tau (\hat{S}_n(t) - S(t))\, d\psi(t)$ for a bounded variation function $\psi$, then the functional delta theorem concludes that $\sqrt{n}(\Psi(\hat{S}_n(t)) - \Psi(S(t)))$ has an asymptotic normal distribution with mean zero and variance $\sigma^2 = \int_0^\tau \int_0^\tau r(s,t)\, d\psi(s)\, d\psi(t)$, where $r(s,t) = E[\mathcal{A}(s;Y,R,L)\mathcal{A}(t;Y,R,L)]$. In this section we want to give a general procedure for estimating $\sigma^2$.

Denote $\mathcal{P}_n$ as the empirical measure of the i.i.d. observations $(Y_i, R_i, L_i)$, $i = 1, \ldots, n$. Clearly, one consistent estimator of $\sigma^2$ is given by

$$\int_0^\tau \int_0^\tau \mathbf{P}_n[\hat{\mathcal{A}}(t;Y,R,L)\hat{\mathcal{A}}(t;Y,R,L)]\, d\psi(s)\, d\psi(t)$$

in which $\hat{\mathcal{A}}(t;Y,R,L)$ is a consistent estimator of $\mathcal{A}(t;Y,R,L)$. To obtain $\hat{\mathcal{A}}(t;Y,R,L)$, we estimate each term in the expression (3.1) separately.

First, in (3.1), we substitute $H_{T|Z^*}(t|Z^*)$ and $H_{C|Z^*}(t|Z^*)$ with their corresponding estimators $\hat{H}_{T|Z^*}(t|\hat{Z})$ and $\hat{H}_{C|Z^*}(t|\hat{Z})$ following Step 3 of Section 2; furthermore, according to the proof of Theorem 3.1, we can consistently estimate the influence functions for $\hat{\beta}_n$ and $\hat{\gamma}_n$, by $\hat{S}(\hat{\beta}_n, Y, R, L)$ and $\hat{S}(\hat{\gamma}_n, Y, R, L)$, respectively. Specifically, $\hat{S}_\beta(\hat{\beta}_n, y, r, l)$ is

$$\left\{\mathcal{P}_n\left[R\left(\frac{\mathbf{P}_n[I_{Y \geq y'}LL'e^{\hat{\gamma}_n'L}]}{\mathbf{P}_n[I_{Y \geq y'}e^{\hat{\gamma}_n'L}]} - \frac{\mathbf{P}_n[I_{Y \geq y'}Le^{\hat{\gamma}_n'L}]^{\otimes 2}}{\mathbf{P}_n[I_{Y \geq y'}e^{\hat{\gamma}_n'L}]^2}\right)\Big|_{y'=Y}\right]\right\}^{-1}$$

$$(3.2) \quad \times \left\{ rl - r \frac{\mathbf{P}_n[I_{y \leq Y} L e^{\hat{\beta}'_n L}]}{\mathbf{P}_n[I_{y \leq Y} e^{\hat{\beta}'_n L}]} - l e^{\hat{\beta}'_n l} \mathbf{P}_n \left[ \frac{R I_{Y \leq y}}{\mathbf{P}_n[I_{Y \leq y'} e^{\hat{\beta}'_n L}]|_{y'=Y}} \right] \right.$$

$$\left. + e^{\hat{\beta}'_n l} \mathbf{P}_n \left[ \frac{R I_{Y \leq y} \mathbf{P}_n[I_{Y \leq y'} L e^{\hat{\beta}'_n L}]|_{y'=Y}}{\mathbf{P}_n[I_{Y \leq y'} e^{\hat{\beta}'_n L}]^2|_{y'=Y}} \right] \right\},$$

and $\hat{S}_\gamma(\hat{\gamma}_n, y, r, l)$ is

$$\left\{ \mathbf{P}_n \left[ (1-R) \left( \frac{\mathbf{P}_n[I_{Y \geq y'} LL' e^{\hat{\gamma}'_n L}]}{\mathbf{P}_n[I_{Y \geq y'} e^{\hat{\gamma}'_n L}]} - \frac{\mathbf{P}_n[I_{Y \geq y'} L e^{\hat{\gamma}'_n L}]^{\otimes 2}}{\mathbf{P}_n[I_{Y \geq y'} e^{\hat{\gamma}'_n L}]^2} \right) \bigg|_{y'=Y} \right] \right\}^{-1}$$

$$(3.3) \quad \times \left\{ (1-r)l - (1-r) \frac{\mathbf{P}_n[I_{y \leq Y} L e^{\hat{\gamma}'_n L}]}{\mathbf{P}_n[I_{y \leq Y} e^{\hat{\gamma}'_n L}]} - l e^{\hat{\gamma}'_n l} \mathbf{P}_n \left[ \frac{(1-R) I_{Y \leq y}}{\mathbf{P}_n[I_{Y \leq y'} e^{\hat{\gamma}'_n L}]|_{y'=Y}} \right] \right.$$

$$\left. + e^{\hat{\gamma}'_n l} \mathbf{P} \left[ \frac{(1-R) I_{Y \leq y} \mathbf{P}_n[I_{Y \leq y'} L e^{\hat{\gamma}'_n L}]|_{y'=Y}}{\mathbf{P}_n[I_{Y \leq y'} e^{\hat{\gamma}'_n L}]^2|_{y'=Y}} \right] \right\}.$$

Additionally, we can estimate

$$(3.4) \quad -E \left[ e^{-H_{T|Z^*}(t|Z^*)} \nabla_\gamma \big|_{\gamma = \gamma^*} \int_0^t \frac{d_u P(T \wedge C \leq u, R = 1 | \gamma' L, \beta^{*'} L)}{P(T \wedge C \geq u | \gamma' L, \beta^{*'} L)} \right]$$

and

$$(3.5) \quad -E \left[ e^{-H_{T|Z^*}(t|Z^*)} \nabla_\beta \big|_{\beta = \beta^*} \int_0^t \frac{d_u P(T \wedge C \leq u, R = 1 | \gamma^{*'} L, \beta' L)}{P(T \wedge C \geq u | \gamma^{*'} L, \beta' L)} \right]$$

using the following lemma.

LEMMA 3.2. *For any constants $(\beta, \gamma)$, we define an estimator of $S(t)$, denoted by $\hat{S}_n(t; \beta, \gamma)$ by repeating Steps 1–4 in Section 2 for fixed $\beta$ and $\gamma$. Let $e_1, \ldots, e_k$ be the canonical bases in $R^{\dim(\beta^*)}$, that is, $e_i$ has 1 at the ith position and 0's elsewhere. Similarly, let $d_1, \ldots, d_l$ be the canonical bases in $R^{\dim(\gamma^*)}$. Moreover, we select a constant $\varepsilon_n$ such that $\varepsilon_n = o(a_n)$, $\sqrt{n} \varepsilon_n \to \infty$. Then when one of the working models is correct, the two statistics, defined by*

$$(3.6) \quad \hat{V}_{\hat{\gamma}_n} = \frac{1}{\varepsilon_n} \begin{pmatrix} \hat{S}_n(t; \hat{\beta}_n, \hat{\gamma}_n + \varepsilon_n d_1) - \hat{S}_n(t) \\ \vdots \\ \hat{S}_n(t; \hat{\beta}_n, \hat{\gamma}_n + \varepsilon_n d_l) - \hat{S}_n(t) \end{pmatrix}$$

*and*

$$(3.7) \quad \hat{V}_{\hat{\beta}_n} = \frac{1}{\varepsilon_n} \begin{pmatrix} \hat{S}_n(t; \hat{\beta}_n + \varepsilon_n e_1, \hat{\gamma}_n) - \hat{S}_n(t) \\ \vdots \\ \hat{S}_n(t; \hat{\beta}_n + \varepsilon_n e_k, \hat{\gamma}_n) - \hat{S}_n(t) \end{pmatrix},$$

*are consistent estimators of* (3.4) *and* (3.5), *respectively.*

So finally, one consistent estimator for $\mathcal{A}(t; Y, R, L)$ is given by

$$e^{-\hat{H}_{T|Z^*}(t|\hat{Z})} - \hat{S}_n(t) - RI_{Y \leq t} e^{\hat{H}_{T|Z^*}(Y|\hat{Z}) + \hat{H}_{C|Z^*}(Y|\hat{Z}) - \hat{H}_{T|Z^*}(t|\hat{Z})}$$
$$+ \int_0^{t \wedge Y} e^{\hat{H}_{T|Z^*}(u|\hat{Z}) + \hat{H}_{C|Z^*}(u|\hat{Z}) - \hat{H}_{T|Z^*}(t|\hat{Z})} \, d_u \hat{H}_{T|Z^*}(u|\hat{Z})$$
$$+ \hat{V}_{\hat{\gamma}_n} \hat{S}_\gamma(\hat{\gamma}_n, Y, R, L) + \hat{V}_{\hat{\beta}_n} \hat{S}_\beta(\hat{\beta}_n, Y, R, L).$$

REMARK 3.3. The numerical method for estimating (3.4) and (3.5) is much more convenient for implementation, compared with the direct estimation of the conditional probabilities in these two expressions. When the bandwidth $a_n$ has order $n^{-1/3}$, one choice of $\varepsilon_n$ may be of the order $n^{-5/12}$. Computationally, except that the final evaluation of the variance requires a numerical double integration, the computing time in the other steps is only a linear order of the computing time for computing $\hat{S}_n(t)$, which is about $O(n^3 a_n^2)$. The storage in the computation is the same order as storing an $n \times n$ numerical array.

REMARK 3.4. As a special example, the asymptotic variance for $\sqrt{n}(\hat{S}_n(t_0) - S(t_0))$ can be approximated by $\mathbf{P}_n[\hat{\mathcal{A}}(t_0; Y, R, L)]^2$ for any fixed time $t_0 \in [0, \tau]$.

**4. Simulation study.** We have performed a simulation study to show the advantages of our approach in small samples. In the simulation, three covariates, denoted as $X_1, X_2, X_3$, were independently generated from the uniform distribution between 0 and 1. The lifetime $T$ was generated from Cox's proportional hazard model, whose hazard rate function had the following form:

$$h_{T|X}(t|X) = \lambda_T(t) \exp\{\beta_1 X_1 + \beta_2 X_2 + \beta_3 X_3 + \beta_{12} X_1 X_2 + \beta_{13} X_1 X_3 + \beta_{23} X_2 X_3\}.$$

The values of the parameters in the simulation were taken to be $\beta_1 = -1$, $\beta_2 = 4$, $\beta_3 = 3$, $\beta_{12} = 0$, $\beta_{13} = 6$, $\beta_{23} = 10$, $\lambda_T(t) = t^4 e^{-5}$. The censoring time $C$ was the minimum of $\tau = 2$ and $C^*$, where $C^*$ was produced using Cox's proportional hazard model with the hazard rate function given by

$$h_{C|X}(t|X) = \lambda_C(t) \exp\{\gamma_1 X_1 + \gamma_2 X_2 + \gamma_3 X_3 + \gamma_{12} X_1 X_2 + \gamma_{13} X_1 X_3 + \gamma_{23} X_2 X_3\}.$$

We chose the parameters as $\gamma_1 = 1$, $\gamma_2 = 1$, $\gamma_3 = 1$, $\gamma_{12} = 0$, $\gamma_{13} = 5$, $\gamma_{23} = 10$, $\lambda_C(t) = t^4 e^{-4.5}$. The choice of the parameter values demonstrated that the dependent censoring between $T$ and $C$ was significant (theoretically, the marginal correlation between $T$ and $C$ was around 75%) and the censoring proportion was not too low (the theoretical censoring probability for this setting is 45%).

We followed the procedure in Section 2 to estimate the survival function for $T$ with the kernel function $k(x_1, x_2) = \exp\{-(x_1^2 + x_2^2)\}$ but started with different working models for $T$ and $C$ given $(X_1, X_2, X_3)$. Especially, if we denoted $\mathbf{X}$ as $(X_1, X_2, X_3)$ and denoted $\mathbf{X}^2$ as their two-way interactions, six pairs of working models could be considered:

Pair 1. We modelled both $T$ and $C$ using all the main effects $\mathbf{X}$ and the two-way interactions $\mathbf{X}^2$ as well as an independent variable $Z$, which was generated from the uniform distribution between 0 and 1.
Pair 2. We modelled both $T$ and $C$ using all the main effects $\mathbf{X}$ and the two-way interactions $\mathbf{X}^2$.
Pair 3. We modelled $T$ using $\mathbf{X}$ and $\mathbf{X}^2$; however, we modelled $C$ using only the main effects $\mathbf{X}$. So we misspecified the model for $C$.
Pair 4. We modelled $C$ using $\mathbf{X}$ and $\mathbf{X}^2$; however, we modelled $T$ using only the main effects $\mathbf{X}$. So we misspecified the model for $T$.
Pair 5. We modelled both $T$ and $C$ using only the main effects $\mathbf{X}$. That is, we misspecified both models.
Pair 6. We did not account for any covariates and the Kaplan–Meier estimate was used to estimate the survival function.

By comparison of the bias and variation among the above six pairs of working models, we expected to verify that the estimates accounting for dependent censoring using covariates in the estimation always perform better than the Kaplan–Meier estimate, that including an irrelevant variable does not bias the estimate and that double robustness is evidenced in small samples.

Moreover, we studied how the estimates varied with the different choices of the sample size $n$, the bandwidth $a_n$ and the oscillation parameter $\varepsilon_n$. We thus generated data with sample size $n = 50$ or $n = 100$. For each generated sample, we used varied bandwidths $a_n = n^{-1/3}, 3n^{-1/3}, 6n^{-1/3}$ to calculate the estimates. In addition, we used different choices $\varepsilon_n = n^{-5/12}, 5n^{-5/12}, 10n^{-5/12}$ to calculate the standard errors and the coverage probabilities in estimating the survival probabilities for $t = \tau/5, 2\tau/5$. Such computation was repeated 500 times.

For both $n = 50$ and $n = 100$, the average censoring proportion was about 45% and the marginal correlation between $T$ and $C$ was 77% in the simulated samples. Tables 1 and 2 report our findings. In Table 1 we give the average mean square error of the estimates on 50 grid points, which is defined as

$$\frac{1}{50}\sum_{i=1}^{50}\left(\hat{S}_n\left(\frac{(i-1)\tau}{50}\right) - S\left(\frac{(i-1)\tau}{50}\right)\right)^2.$$

In Table 2 we report the average bias and the 95% confidence interval coverage probabilities for estimating the survival probabilities at times $\tau/5$ and $2\tau/5$. Since it has been found that the coverage probabilities vary very little when $\varepsilon_n$ varies in our choices, we only report the results for $\varepsilon_n = n^{-5/12}$.

From Table 1, it is clear that the Kaplan–Meier estimates have the largest mean square error and the estimates adjusting for dependent censoring using covariates can reduce it by 50% for sample size 50 and by over 60% for sample size 100. Moreover, using the irrelevant covariate $Z$ in the regression models does not increase the mean square error, and when either of the regression models is correct (i.e., both the main effects and the two-way interactions among $X_1, X_2, X_3$ are used in the regression), the mean square errors are, on average, 10% less than for the case which only uses the main effects in both regressions. The mean square errors of the estimates are fairly robust to the choice of the bandwidth. The results displayed in Table 2 further evidence the above findings from the view of the point estimates of $S(t)$ and the corresponding coverage probabilities. Table 2 shows that when either regression model is specified correctly, the biases in the estimates are less than for the cases when both models are misspecified; the Kaplan–Meier induces the largest biases. Overall, these biases decrease by 50% when the sample size increases from 50 to 100. With the sample size 50 or 100, the coverage probabilities using the methods proposed in Section 3 are fairly accurate for $t = \tau/5$ when either regression model is specified correctly; however, they tend to be smaller for $t = 2\tau/5$ due to the larger bias caused by high censoring at the tail. When the bandwidth is large (for instance, $a_n = 6n^{-1/3}$), the biases increase due to oversmoothing, but the coverage probabilities do not vary much.

Our simulation study indicates that the estimates of the survival function by adjusting for dependent censoring using auxiliary covariates always

TABLE 1
*Mean square error from 500 samples*

| $n$ | model $T$ | model $C$ | MSE($\times 10^{-3}$) $a_n = n^{-1/3}$ | MSE($\times 10^{-3}$) $a_n = 3n^{-1/3}$ | MSE($\times 10^{-3}$) $a_n = 6n^{-1/3}$ |
|---|---|---|---|---|---|
| 50 | $(\mathbf{X}, \mathbf{X}^2, Z)$ | $(\mathbf{X}, \mathbf{X}^2, Z)$ | 7.2 | 6.9 | 6.8 |
|  | $(\mathbf{X}, \mathbf{X}^2)$ | $(\mathbf{X}, \mathbf{X}^2)$ | 7.2 | 6.8 | 6.8 |
|  | $(\mathbf{X}, \mathbf{X}^2)$ | $(\mathbf{X})^{\text{a}}$ | 7.0 | 6.7 | 6.7 |
|  | $(\mathbf{X})^{\text{a}}$ | $(\mathbf{X}, \mathbf{X}^2)$ | 7.1 | 6.9 | 7.0 |
|  | $(\mathbf{X})^{\text{a}}$ | $(\mathbf{X})^{\text{a}}$ | 7.7 | 7.4 | 7.5 |
|  | $(-)^{\text{b}}$ | $(-)^{\text{b}}$ | 17.4 | 17.4 | 17.4 |
| 100 | $(\mathbf{X}, \mathbf{X}^2, Z)$ | $(\mathbf{X}, \mathbf{X}^2, Z)$ | 3.5 | 3.4 | 3.3 |
|  | $(\mathbf{X}, \mathbf{X}^2)$ | $(\mathbf{X}, \mathbf{X}^2)$ | 3.5 | 3.4 | 3.3 |
|  | $(\mathbf{X}, \mathbf{X}^2)$ | $(\mathbf{X})^{\text{a}}$ | 3.4 | 3.3 | 3.3 |
|  | $(\mathbf{X})^{\text{a}}$ | $(\mathbf{X}, \mathbf{X}^2)$ | 3.7 | 3.5 | 3.5 |
|  | $(\mathbf{X})^{\text{a}}$ | $(\mathbf{X})^{\text{a}}$ | 4.3 | 4.0 | 4.0 |
|  | $(-)^{\text{b}}$ | $(-)^{\text{b}}$ | 13.0 | 13.0 | 13.0 |

Notation. $(\cdots)^{\text{a}}$: model is misspecified; $(-)^{\text{b}}$: the Kaplan–Meier estimate is used.

TABLE 2
*Estimate of the survival probability at times $t = \tau/5$ and $t = 2\tau/5$ from 500 samples with $\varepsilon_n = n^{-5/12}$*

| $n$ | $a_n$ | model $T$ | model $C$ | $S(\tau/5)$ bias($\times 10^{-2}$) | 95% cp | $S(2\tau/5)$ bias($\times 10^{-2}$) | 95% cp |
|---|---|---|---|---|---|---|---|
| 50 | $n^{-1/3}$ | $(\mathbf{X}, \mathbf{X}^2, Z)$ | $(\mathbf{X}, \mathbf{X}^2, Z)$ | 0.94 | 0.94 | 1.97 | 0.92 |
| | | $(\mathbf{X}, \mathbf{X}^2)$ | $(\mathbf{X}, \mathbf{X}^2)$ | 0.94 | 0.94 | 1.95 | 0.92 |
| | | $(\mathbf{X}, \mathbf{X}^2)$ | $(\mathbf{X})^{\text{a}}$ | 0.96 | 0.94 | 2.03 | 0.93 |
| | | $(\mathbf{X})^{\text{a}}$ | $(\mathbf{X}, \mathbf{X}^2)$ | 1.05 | 0.93 | 2.23 | 0.92 |
| | | $(\mathbf{X})^{\text{a}}$ | $(\mathbf{X})^{\text{a}}$ | 1.22 | 0.93 | 2.99 | 0.92 |
| | | $(-)^{\text{b}}$ | $(-)^{\text{b}}$ | 5.69 | 0.87 | 10.72 | 0.70 |
| | $3n^{-1/3}$ | $(\mathbf{X}, \mathbf{X}^2, Z)$ | $(\mathbf{X}, \mathbf{X}^2, Z)$ | 0.94 | 0.94 | 1.95 | 0.91 |
| | | $(\mathbf{X}, \mathbf{X}^2)$ | $(\mathbf{X}, \mathbf{X}^2)$ | 0.96 | 0.94 | 2.02 | 0.91 |
| | | $(\mathbf{X}, \mathbf{X}^2)$ | $(\mathbf{X})^{\text{a}}$ | 1.00 | 0.93 | 2.12 | 0.90 |
| | | $(\mathbf{X})^{\text{a}}$ | $(\mathbf{X}, \mathbf{X}^2)$ | 1.08 | 0.92 | 2.23 | 0.89 |
| | | $(\mathbf{X})^{\text{a}}$ | $(\mathbf{X})^{\text{a}}$ | 1.46 | 0.92 | 3.02 | 0.87 |
| | | $(-)^{\text{b}}$ | $(-)^{\text{b}}$ | 5.69 | 0.87 | 10.72 | 0.70 |
| | $6n^{-1/3}$ | $(\mathbf{X}, \mathbf{X}^2, Z)$ | $(\mathbf{X}, \mathbf{X}^2, Z)$ | 1.30 | 0.93 | 2.50 | 0.90 |
| | | $(\mathbf{X}, \mathbf{X}^2)$ | $(\mathbf{X}, \mathbf{X}^2)$ | 1.33 | 0.93 | 2.58 | 0.90 |
| | | $(\mathbf{X}, \mathbf{X}^2)$ | $(\mathbf{X})^{\text{a}}$ | 1.50 | 0.93 | 2.72 | 0.89 |
| | | $(\mathbf{X})^{\text{a}}$ | $(\mathbf{X}, \mathbf{X}^2)$ | 1.54 | 0.92 | 2.76 | 0.89 |
| | | $(\mathbf{X})^{\text{a}}$ | $(\mathbf{X})^{\text{a}}$ | 2.16 | 0.92 | 3.58 | 0.87 |
| | | $(-)^{\text{b}}$ | $(-)^{\text{b}}$ | 5.69 | 0.87 | 10.72 | 0.70 |
| 100 | $n^{-1/3}$ | $(\mathbf{X}, \mathbf{X}^2, Z)$ | $(\mathbf{X}, \mathbf{X}^2, Z)$ | 0.44 | 0.94 | 1.37 | 0.96 |
| | | $(\mathbf{X}, \mathbf{X}^2)$ | $(\mathbf{X}, \mathbf{X}^2)$ | 0.44 | 0.94 | 1.33 | 0.95 |
| | | $(\mathbf{X}, \mathbf{X}^2)$ | $(\mathbf{X})^{\text{a}}$ | 0.48 | 0.93 | 1.36 | 0.94 |
| | | $(\mathbf{X})^{\text{a}}$ | $(\mathbf{X}, \mathbf{X}^2)$ | 0.52 | 0.92 | 1.51 | 0.92 |
| | | $(\mathbf{X})^{\text{a}}$ | $(\mathbf{X})^{\text{a}}$ | 0.95 | 0.93 | 2.73 | 0.89 |
| | | $(-)^{\text{b}}$ | $(-)^{\text{b}}$ | 5.44 | 0.75 | 10.50 | 0.51 |
| | $3n^{-1/3}$ | $(\mathbf{X}, \mathbf{X}^2, Z)$ | $(\mathbf{X}, \mathbf{X}^2, Z)$ | 0.45 | 0.92 | 1.46 | 0.93 |
| | | $(\mathbf{X}, \mathbf{X}^2)$ | $(\mathbf{X}, \mathbf{X}^2)$ | 0.48 | 0.93 | 1.42 | 0.92 |
| | | $(\mathbf{X}, \mathbf{X}^2)$ | $(\mathbf{X})^{\text{a}}$ | 0.57 | 0.92 | 1.48 | 0.92 |
| | | $(\mathbf{X})^{\text{a}}$ | $(\mathbf{X}, \mathbf{X}^2)$ | 0.55 | 0.91 | 1.52 | 0.91 |
| | | $(\mathbf{X})^{\text{a}}$ | $(\mathbf{X})^{\text{a}}$ | 1.13 | 0.90 | 2.83 | 0.89 |
| | | $(-)^{\text{b}}$ | $(-)^{\text{b}}$ | 5.44 | 0.75 | 10.50 | 0.51 |
| | $6n^{-1/3}$ | $(\mathbf{X}, \mathbf{X}^2, Z)$ | $(\mathbf{X}, \mathbf{X}^2, Z)$ | 0.79 | 0.92 | 1.89 | 0.91 |
| | | $(\mathbf{X}, \mathbf{X}^2)$ | $(\mathbf{X}, \mathbf{X}^2)$ | 0.82 | 0.92 | 1.89 | 0.92 |
| | | $(\mathbf{X}, \mathbf{X}^2)$ | $(\mathbf{X})^{\text{a}}$ | 0.99 | 0.93 | 2.09 | 0.91 |
| | | $(\mathbf{X})^{\text{a}}$ | $(\mathbf{X}, \mathbf{X}^2)$ | 0.96 | 0.91 | 1.98 | 0.91 |
| | | $(\mathbf{X})^{\text{a}}$ | $(\mathbf{X})^{\text{a}}$ | 1.67 | 0.90 | 3.23 | 0.89 |
| | | $(-)^{\text{b}}$ | $(-)^{\text{b}}$ | 5.44 | 0.75 | 10.50 | 0.51 |

Notation. $(\cdots)^{\text{a}}$: model is misspecified; $(-)^{\text{b}}$: the Kaplan–Meier estimate is used.

induce smaller mean square errors, fewer biases and more accurate coverage probabilities compared with the Kaplan–Meier estimates. Moreover, the estimates have better performance when either the model for $T$ or the model for $C$ given the covariates is used correctly. The overall mean square errors of the consistent estimates are fairly robust to the choice of the bandwidth; but the point estimates and the inference vary with the choices of the bandwidth and the location of time points.

**5. Discussion.** Both our theoretical justification of large samples and simulation studies with small samples conclude that, when right-censored data include high-dimensional auxiliary covariates, condensing such information by utilizing working models for both lifetime and censoring time given covariates can make adjusting for dependent censoring possible and produce an estimator which is robust to the misspecification of either working model and robust to accidentally using irrelevant information.

It is observed in our simulations that the choice of the bandwidth $a_n$ plays an important role in influencing the bias and the inference for the point estimate. A large $a_n$ may oversmooth the conditional hazard rate estimator (in fact, with simple calculation, for fixed $n$, if $a_n$ is close to infinity, our estimate approximates the Kaplan–Meier estimate), while a small $a_n$ may overfit the conditional hazard rate estimator, and thus introduce large variation in estimation. So far, we let $a_n$ be a constant only depending on $n$ and no general selection rule is followed; however, the simulation results imply that a data-adaptive and location-adaptive $a_n$ may give a better performing estimate. The cross-validation approach may be used to choose $a_n$ or we can use the $k$ nearest neighbor approach in nonparametric hazard regression. We will explore this issue more in the future.

Though we hope that our working models are correct, we never know in reality. To make this hope more likely, we may use more general models other than Cox proportional hazard models as working models, for example, we can use a generalized additive model or use splines as covariates in the working models, and so on. A model selection rule is thus useful in choosing the optimal working models in terms of the performance of the estimates and the model complexity. Therefore, a test for goodness of fit as well as a test for comparing two different sets of working models will be useful in practice.

Finally, when $L$ includes the time-dependent covariates, our approach is not obvious to fit this situation. This is because the condensed information $(\beta'L, \gamma'L)$ is still time-dependent so their dimension is infinite; then an essential problem is how to derive a nonparametric estimate of the marginal survival function in the presence of even a single time-dependent covariate. Further exploration of this issue is ongoing.

## APPENDIX: PROOFS

PROOF OF THEOREM 3.1. We consider only the estimator $\hat{\beta}_n$ in the following. The argument for the estimator $\hat{\gamma}_n$ is similar. Obviously, $\hat{\beta}_n$ maximizes

$$\widetilde{L}_1^{(n)}(\beta) = \frac{1}{n}\sum_{i=1}^n R_i \beta' L_i - \frac{1}{n}\sum_{i=1}^n R_i \log\left(\sum_{Y_j \geq Y_i} e^{\beta' L_j}\right).$$

Note that $\widetilde{L}_1^n(\beta)$ is a concave function of $\beta$ and its limit, which is equal to $\widetilde{L}_1(\beta) = \mathbf{P}[R\beta'L - R\log E[I_{Y \geq y}e^{\beta'L}]|_{y=Y}]$, is a strictly concave function. By an argument similar to that in Andersen and Gill (1982), we obtain that with probability 1, $\hat{\beta}_n$ converges to the unique maximum of $\widetilde{L}_1(\beta)$, denoted by $\beta^*$.

After the linearization of the equation $\tilde{L}_1^{(n)}(\hat{\beta}_n) = 0$ around $\beta^*$, we obtain that

$$\sqrt{n}(\hat{\beta}_n - \beta^*) = \sqrt{n}(\mathbf{P}_n - \mathbf{P})S(\beta^*, Y, R, L) + o_p(1),$$

where the influence function $S(\beta^*, y, r, l)$ is equal to

$$-\{\nabla_{\beta\beta}^2 \widetilde{L}_1(\beta^*)\}^{-1}\left\{rl - r\frac{\mathbf{P}[I_{y \leq Y}Le^{\beta^{*\prime}L}]}{\mathbf{P}[I_{y \leq Y}e^{\beta^{*\prime}L}]} - le^{\beta^{*\prime}l}\mathbf{P}\left[\frac{RI_{Y \leq y}}{\mathbf{P}[I_{Y \leq y'}e^{\beta^{*\prime}L}]|_{y'=Y}}\right]\right.$$
$$\left. + e^{\beta^{*\prime}l}\mathbf{P}\left[\frac{RI_{Y \leq y}\mathbf{P}[I_{Y \leq y'}Le^{\beta^{*\prime}L}]|_{y'=Y}}{\mathbf{P}[I_{Y \leq y'}e^{\beta^{*\prime}L}]^2|_{y'=Y}}\right]\right\}. \quad \square$$

PROOF OF THEOREM 3.2. Recall $Z^* = (\beta^{*\prime}L, \gamma^{*\prime}L)$ and $\hat{Z} = (\hat{\beta}_n'L, \hat{\gamma}_n'L)$. We assume one of the working models is correct so $T$ and $C$ are independent given $Z^*$ from Lemma 3.1. The whole proof consists of three steps: In the first step, we show the uniform consistency of $\hat{H}_{T|Z^*}(t|z)$, thus $\hat{S}_n(t)$; then we write $\sqrt{n}(\hat{S}_n(t) - S(t))$ as a linear functional of the empirical processes; in the third step, we apply empirical process theory to obtain the asymptotic properties.

First, the following result holds and the proof is given in Dabrowska (1987).

LEMMA A.1. *For any $z$ in the support of $Z^*$,*

$$\left\|\frac{\mathbf{P}_n[K((Z^* - z)/a_n)I_{Y \geq t}]}{\mathbf{P}_n[K((Z^* - z)/a_n)]} - P(Y \geq t|Z^* = z)\right\|_{l^\infty([0,\tau])} \xrightarrow{P} 0,$$

$$\left\|\frac{\mathbf{P}_n[K((Z^* - z)/a_n)I_{Y \leq t}R]}{\mathbf{P}_n[K((Z^* - z)/a_n)]} - P(Y \leq t, R = 1|Z^* = z)\right\|_{l^\infty([0,\tau])} \xrightarrow{P} 0.$$

LEMMA A.2. *For any $z$ in the support of $Z^*$,*

$$\left\|\frac{\mathbf{P}_n[K((\hat{Z}-z)/a_n)I_{Y\geq t}]}{\mathbf{P}_n[K((\hat{Z}-z)/a_n)]} - P(Y\geq t|Z^*=z)\right\|_{l^\infty([0,\tau])} \xrightarrow{P} 0,$$

$$\left\|\frac{\mathbf{P}_n[K((\hat{Z}-z)/a_n)I_{Y\leq t}R]}{\mathbf{P}_n[K((\hat{Z}-z)/a_n)]} - P(Y\leq t, R=1|Z^*=z)\right\|_{l^\infty([0,\tau])} \xrightarrow{P} 0.$$

PROOF. For convenience, we denote

$$g_n(\beta,\gamma) = \frac{1}{a_n^2}\mathbf{P}_n\left[K\left(\frac{(\beta'L,\gamma'L)-z}{a_n}\right)I_{Y\geq t}\right].$$

We show that $\sup_{t\in[0,\tau]}|g_n(\hat{\beta}_n,\hat{\gamma}_n) - g_n(\beta^*,\gamma^*)| \to 0$ a.s. By the property of the kernel function and the mean value theorem, we have that

$$|g_n(\hat{\beta}_n,\hat{\gamma}_n) - g_n(\beta^*,\gamma^*)|$$
$$\leq \frac{1}{na_n^2}\sum_{i=1}^n\left|\nabla K\left(\frac{(\tilde{\beta}'L,\tilde{\gamma}'L)-z}{a_n}\right)\right|O\left(\frac{|\hat{\beta}_n-\beta^*|}{a_n}+\frac{|\hat{\gamma}_n-\gamma^*|}{a_n}\right),$$

where $(\tilde{\beta},\tilde{\gamma})$ is between $(\hat{\beta}_n,\hat{\gamma}_n)$ and $(\beta^*,\gamma^*)$. Hence, for any $z=(z_1,z_2)$ in the support of $Z^*$,

$$|g_n(\hat{\beta}_n,\hat{\gamma}_n) - g_n(\beta^*,\gamma^*)|$$
$$\leq O_p\left(\frac{1}{\sqrt{n}a_n}\right)\left[\frac{1}{na_n^2}\sum_{i=1}^n\frac{1}{1+(\tilde{\beta}'L_i-z_1)^2/a_n^2+(\tilde{\gamma}'L_i-z_2)^2/a_n^2}\right]$$
$$\leq O_p\left(\frac{1}{\sqrt{n}a_n}\right)\left[\frac{1}{na_n^2}\sum_{i=1}^n\frac{1}{1+(\beta^{*'}L_i-z_1)^2/a_n^2+(\gamma^{*'}L_i-z_2)^2/a_n^2}\right],$$

where the last step follows because $|\nabla_{x_j}\log(1+x_1^2+x_2^2)|, j=1,2$, is uniformly bounded and $\frac{|\tilde{\beta}'L-\beta^{*'}L|}{a_n}+\frac{|\tilde{\gamma}'L-\gamma^{*'}L|}{a_n} \leq O_p(1)$. Notice that

$$\mathbf{P}\left[\frac{1}{a_n^2}\frac{1}{1+(\beta^{*'}L-z_1)^2/a_n^2+(\gamma^{*'}L-z_2)^2/a_n^2}\right]$$

is uniformly bounded. So $\sup_{t\in[0,\tau]}|g_n(\hat{\beta}_n,\hat{\gamma}_n) - g_n(\beta^*,\gamma^*)| \leq O_p(\frac{1}{\sqrt{n}a_n})$.

Similarly, we can obtain that

$$\sup_{t\in[0,\tau]}\left|\frac{1}{na_n^2}\sum_{i=1}^n K\left(\frac{\hat{Z}-z}{a_n}\right) - \frac{1}{na_n^2}\sum_{i=1}^n K\left(\frac{Z^*-z}{a_n}\right)\right| \xrightarrow{p} 0.$$

Combining this result with Lemma A.1, it is clear the first half of Lemma A.2 holds. The second half of Lemma A.2 can be proved similarly. □

LEMMA A.3. *Denote*

$$\hat{H}_{T|Z^*}(t|z) = \int_0^t \frac{d_s \mathbf{P}_n[K((\hat{Z}-z)/a_n)RI_{Y \leq s}]}{\mathbf{P}_n[K((\hat{Z}-z)/a_n)I_{Y \geq s}]}$$

*and*

$$\hat{S}_{T|Z^*}(t|z) = \prod_{s \leq t}(1 - \hat{H}_{T|Z^*}(\{s\}|z)).$$

*Then for any $z$ in the support of $Z^*$, in probability $\|\hat{H}_{T|Z^*}(t|z) - H_{T|Z^*}(t|z)\|_{l^\infty([0,\tau])} \to 0$ and $\|\hat{S}_{T|Z^*}(t|z) - S_{T|Z^*}(t|z)\|_{l^\infty([0,\tau])} \to 0$.*

PROOF. The first result follows from Assumption 3.3, Lemma A.2 and the following inequality:

$$\|\hat{H}_{T|Z^*}(t|z) - H_{T|Z^*}(t|z)\|_{l^\infty([0,\tau])}$$

$$= \left\| \int_0^t \frac{d_s \mathbf{P}_n[K((\hat{Z}-z)/a_n)RI_{Y \leq s}]}{\mathbf{P}_n[K((\hat{Z}-z)/a_n)I_{Y \geq s}]} - \int_0^t \frac{d_s E[RI_{Y \leq s}|Z^*]}{E[I_{Y \geq s}|Z^*]} \right\|_{l^\infty([0,\tau])}$$

$$\leq \left\| \frac{\mathbf{P}_n[K((\hat{Z}-z)/a_n)I_{Y \leq t}R]}{\mathbf{P}_n[K((\hat{Z}-z)/a_n)]} - P(Y \leq t, R = 1|Z^* = z) \right\|_{l^\infty([0,\tau])}$$

$$\times \left\{ \min_{t \in [0,\tau]} \left| \frac{\mathbf{P}_n[K((\hat{Z}-z)/a_n)I_{Y \geq t}]}{\mathbf{P}_n[K((\hat{Z}-z)/a_n)]} \right| \right\}^{-1}$$

$$+ \left\| \frac{\mathbf{P}_n[K((\hat{Z}-z)/a_n)I_{Y \geq t}]}{\mathbf{P}_n[K((\hat{Z}-z)/a_n)]} - P(Y \geq t|Z^* = z) \right\|_{l^\infty([0,\tau])}$$

$$\times \left\{ \min_{t \in [0,\tau]} \left| \frac{\mathbf{P}_n[K((\hat{Z}-z)/a_n)I_{Y \geq t}]}{\mathbf{P}_n[K((\hat{Z}-z)/a_n)]} P(Y \geq t|Z^* = z) \right| \right\}^{-1}. \quad \square$$

For the second result, we use the Duhamel equation and integration by parts: for any $t \in [0, \tau]$,

$$|\hat{S}_{T|Z^*}(t|z) - S_{T|Z^*}(t|z)|$$

$$= \left| S_{T|Z^*}(t|z) \int_0^t \frac{\hat{S}_{T|Z^*}(u-|z)}{S_{T|Z^*}(u|z)} d(\hat{H}_{T|Z^*}(u|z) - H_{T|Z^*}(u|z)) \right|$$

$$= \left| \hat{S}_{T|Z^*}(t-|z)(\hat{H}_{T|Z^*}(t|z) - H_{T|Z^*}(t|z)) \right.$$

$$\left. - \int_0^t (\hat{H}_{T|Z^*}(u|z) - H_{T|Z^*}(u|z)) \right.$$

MARGINAL SURVIVAL ESTIMATION 19$$\times \left( \frac{d\hat{S}_{T|Z^*}(u-|z)}{S_{T|Z^*}(u|z)} - \frac{\hat{S}_{T|Z^*}(u-|z)\, dS_{T|Z^*}(u|z)}{S_{T|Z^*}(u|z)^2} \right) \Bigg|$$

$$\leq \left( 1 + \frac{2}{\min P(T > \tau | Z^* = z)^2} \right) \max_{0 \leq s \leq t} |\hat{H}_{T|Z^*}(s|z) - H_{T|Z^*}(s|z)|.$$

For the second step, we will write $\hat{H}_{T|Z^*}(t|z) - H_{T|Z^*}(t|z)$, then $\hat{S}_n(t) - S(t)$ in terms of the empirical process $(\mathbf{P}_n - \mathbf{P})$. First, we obtain

$$\hat{H}_{T|Z^*}(t|z) - H_{T|Z^*}(t|z)$$

$$= (\mathbf{P}_n - \mathbf{P})\left[ \frac{1/a_n^2 K((\hat{Z}-z)/a_n) I_{Y \leq t} R}{\mathbf{P}_n[1/a_n^2 K((\hat{Z}-z)/a_n) I_{Y \geq y}]|_{y=Y}} \right]$$

$$- \mathbf{P}\left[ \frac{1/a_n^2 K((\hat{Z}-z)/a_n) I_{Y \leq t} R (\mathbf{P}_n - \mathbf{P})[1/a_n^2 K((\hat{Z}-z)/a_n) I_{Y \geq y}]|_{y=Y}}{\mathbf{P}_n[1/a_n^2 K((\hat{Z}-z)/a_n) I_{Y \geq y}]^2|_{y=Y}} \right]$$

$$+ \left\{ \mathbf{P}\left[ \frac{1/a_n^2 K((\hat{Z}-z)/a_n) I_{Y \leq t} R}{\mathbf{P}[1/a_n^2 K((\hat{Z}-z)/a_n) I_{Y \geq y}]|_{y=Y}} \right] - H_{T|Z^*}(t|z) \right\}$$

$$= I + II + III.$$

For $III$, a simple transformation in the integral gives that uniformly in $z$ in the support of $Z^*$ and $t \in [0, \tau]$,

$$III = \int_0^t \frac{d_u P(R=1, Y \leq u | \hat{Z} = z)}{P(Y \geq u | \hat{Z} = z)} + O_p(a_n^2) - H_{T|Z^*}(t|z).$$

On the other hand, since by Lemma 3.1 $\int_0^t \frac{d_u P(R=1, Y \leq u | Z^* = z)}{P(Y \geq u | Z^* = z)} = H_{T|Z^*}(t|z)$, we perform the Taylor expansion of the above expansion around $(\beta^*, \gamma^*)$, and then $III$ becomes

$$III = \nabla_\beta|_{\beta=\beta^*}\left[ \int_0^t \frac{d_u P(Y \leq u, R=1 | (\beta'L, \gamma^{*'}L) = z)}{P(Y \geq u | (\beta'L, \gamma^{*'}L) = z)} \right](\hat{\beta}_n - \beta^*)$$

$$+ \nabla_\gamma|_{\gamma=\gamma^*}\left[ \int_0^t \frac{d_u P(Y \leq u, R=1 | (\beta^{*'}L, \gamma'L) = z)}{P(Y \geq u | (\beta^{*'}L, \gamma'L) = z)} \right](\hat{\gamma}_n - \gamma^*)$$

$$+ O_p(a_n^2) + O_p\left(\frac{1}{n}\right).$$

For convenience, we introduce more notation:

$$h_1^n(y, r, l; \beta, \gamma, t, z) = \frac{1/a_n^2 K(((\beta'l, \gamma'l) - z)/a_n) I_{y \leq t} r}{\mathbf{P}_n[1/a_n^2 K(((\beta'L, \gamma'L) - z)/a_n) I_{Y \geq y}]},$$

$$h_2(y, l; \beta, \gamma, t, z) = \frac{1}{a_n^2} K\left( \frac{(\beta'l, \gamma'l) - z}{a_n} \right)$$

$$\times \mathbf{P}\left[\frac{1/a_n^2 K(((\beta'L, \gamma'L) - z)/a_n) I_{Y \leq t} I_{Y \leq y} R}{\mathbf{P}_n[1/a_n^2 K(((\beta'L, \gamma'L) - z)/a_n) I_{Y \geq y}]^2|_{y=Y}}\right],$$

$$B(\beta, \gamma, z, t) = \int_0^t \frac{d_u P(Y \leq u, R = 1|(\beta'L, \gamma'L) = z)}{P(Y \geq u|(\beta'L, \gamma'L) = z)}.$$

After substituting this notation into the expression $\hat{H}_{T|Z^*}(t|z) - H_{T|Z^*}(t|z)$, then further substituting into $\hat{S}_{T|Z^*}(t|z) - S_{T|Z^*}(t|z)$ in the Duhamel equation, we have that, uniformly in $t \in [0, \tau]$,

$$\hat{S}_{T|Z^*}(t|z) - S_{T|Z^*}(t|z)$$

$$= -S(t|z)\bigg\{(\mathbf{P}_n - \mathbf{P})\bigg[\int_0^t \frac{\hat{S}_{T|Z^*}(u-|z)}{S_{T|Z^*}(u|z)} dh_1^n(Y, R, L; \hat{\beta}_n, \hat{\gamma}_n, u, z)\bigg]$$

$$- (\mathbf{P}_n - \mathbf{P})\bigg[\int_0^t \frac{\hat{S}_{T|Z^*}(u-|z)}{S_{T|Z^*}(u|z)} dh_2^n(Y, L; \hat{\beta}_n, \hat{\gamma}_n, u, z)\bigg]$$

(A.1)

$$+ \bigg[\int_0^t \frac{\hat{S}_{T|Z^*}(u-|z)}{S_{T|Z^*}(u|z)} d_u \nabla_\beta B(\beta^*, \gamma^*, z, u)\bigg](\hat{\beta}_n - \beta^*)$$

$$+ \bigg[\int_0^t \frac{\hat{S}_{T|Z^*}(u-|z)}{S_{T|Z^*}(u|z)} d_u \nabla_\gamma B(\beta^*, \gamma^*, z, u)\bigg](\hat{\gamma}_n - \gamma^*)\bigg\}$$

$$+ O_p(a_n^2) + O_p\bigg(\frac{1}{n}\bigg).$$

Note that

$$\sqrt{n}(\hat{S}_n(t) - S(t)) = \sqrt{n}(\mathbf{P}_n[\hat{S}_{T|Z^*}(t|\hat{Z})] - \mathbf{P}[S_{T|Z^*}(t|\hat{Z})])$$

$$= \sqrt{n}(\mathbf{P}_n - \mathbf{P})[\hat{S}_{T|Z^*}(t|\hat{Z})] + \sqrt{n}\mathbf{P}[(\hat{S}_{T|Z^*}(t|\hat{Z}) - S(t|\hat{Z}))].$$

After using (A.1) and the results of Lemmas A.2 and A.3, we obtain that uniformly in $t \in [0, \tau]$,

(A.2)
$$\sqrt{n}(\hat{S}_n(t) - S(t))$$
$$= \sqrt{n}(\mathbf{P}_n - \mathbf{P})w_n(Y, R, L; \hat{\beta}_n, \hat{\gamma}_n, t)$$
$$+ \mathbf{P}[S_{T|Z^*}(t|Z^*)\nabla_\beta B(\beta^*, \gamma^*, Z^*, t)]\sqrt{n}(\hat{\beta}_n - \beta^*)$$
$$+ \mathbf{P}[S_{T|Z^*}(t|Z^*)\nabla_\gamma B(\beta^*, \gamma^*, Z^*, t)]\sqrt{n}(\hat{\gamma}_n - \gamma^*) + o_p(1),$$

where

$$w_n(y, r, l; \hat{\beta}_n, \hat{\gamma}_n, t)$$
$$= \hat{S}_{T|Z^*}(t|\hat{z}) - S(t)$$

$$- \mathbf{P}\left[S_{T|Z^*}(t|\hat{Z}) \int_0^t \frac{\hat{S}_{T|Z^*}(u-|\hat{Z})}{S_{T|Z^*}(u|\hat{Z})} \, dh_1^n(y,r,l;\hat{\beta}_n,\hat{\gamma}_n,u,\hat{Z})\right]$$

$$+ \mathbf{P}\left[S_{T|Z^*}(t|\hat{Z}) \int_0^t \frac{\hat{S}_{T|Z^*}(u-|\hat{Z})}{S_{T|Z^*}(u|\hat{Z})} \, dh_2^n(y,l;\hat{\beta}_n,\hat{\gamma}_n,u,\hat{Z})\right].$$

In the third step, empirical process theory is applied to the above expression for $\sqrt{n}(\hat{S}_n(t) - S(t))$ to obtain the asymptotic properties of $\hat{S}_n(t)$. We consider the empirical process

$$\left\{\sqrt{n}(\mathbf{P}_n - \mathbf{P})w_n\left(Y, R, L; \beta^* + \frac{\theta_1}{\sqrt{n}}, \gamma^* + \frac{\theta_2}{\sqrt{n}}, t\right) : t \in [0, \tau],\right.$$
$$\left.\theta_1 = O_p(1), \theta_2 = O_p(1)\right\},$$

which is indexed by $(t, \theta_1, \theta_2)$. First, we claim that uniformly in $t$,

$$w_n(Y, R, L; \hat{\beta}_n, \hat{\gamma}_n, t)$$
$$\to S_{T|Z^*}(t|Z^*) - S(t) - \frac{RI_{Y \leq t}S_{T|Z^*}(t|Z^*)}{P(Y \geq y'|Z^*)|_{y'=Y}}$$
$$+ S_{T|Z^*}(t|Z^*) \int_0^{t \wedge Y} e^{H_{T|Z^*}(u|Z^*) + H_{C|Z^*}(u|Z^*)} \, dH_{T|Z^*}(u|Z^*)$$

in probability. This is true by using arguments similar to those in the proofs of Lemmas A.2 and A.3. Second, with technical calculation we can verify that each function in the class indexed by $(t, \theta_1, \theta_2)$ belongs to $BV[0, \tau]$ as a function of $t$ and is Lipschitz continuous with respect to $(\theta_1, \theta_2)$ with the Lipschitz coefficient bounded by $O(\frac{1}{\sqrt{na_n}})$ in probability. Thus we can check each condition of Theorem 2.11.23 in van der Vaart and Wellner (1996) and obtain that, in $l^\infty([0, \tau])$,

$$\sqrt{n}(\mathbf{P}_n - \mathbf{P})w_n(Y, R, L; \hat{\beta}_n, \hat{\gamma}_n, t)$$
$$= \sqrt{n}(\mathbf{P}_n - \mathbf{P})$$
$$\times \left[S_{T|Z^*}(t|Z^*) - S(t) - \frac{RI_{Y \leq t}S_{T|Z^*}(t|Z^*)}{P(Y \geq y'|Z^*)|_{y'=Y}}\right.$$
$$\left. + S_{T|Z^*}(t|Z^*) \int_0^{t \wedge Y} e^{H_{T|Z^*}(u|Z^*) + H_{C|Z^*}(u|Z^*)} \, d_u H_{T|Z^*}(u|Z^*)\right]$$
$$+ o_p(1).$$

Therefore, from (A.2) we obtain that uniformly in $t \in [0, \tau]$,

$$\sqrt{n}(\hat{S}_n(t) - S(t))$$

22    D. ZENG$$= \sqrt{n}(\mathbf{P}_n - \mathbf{P})$$
$$\times \left[ S_{T|Z^*}(t|Z^*) - S(t) - \frac{RI_{Y\leq t}S_{T|Z^*}(t|Z^*)}{P(Y\geq y'|Z^*)|_{y'=Y}} \right.$$
$$\left. + S_{T|Z^*}(t|Z^*) \int_0^{t\wedge Y} e^{H_{T|Z^*}(u|Z^*)+H_{C|Z^*}(u|Z^*)} d_u H_{T|Z^*}(u|Z^*) \right]$$
$$- \mathbf{P}[S_{T|Z^*}(t|Z^*)\nabla_\beta B(\beta^*,\gamma^*,Z^*,t)]\sqrt{n}(\hat\beta_n - \beta^*)$$
$$- \mathbf{P}[S_{T|Z^*}(t|Z^*)\nabla_\gamma B(\beta^*,\gamma^*,Z^*,t)]\sqrt{n}(\hat\gamma_n - \gamma^*) + o_p(1).$$

Combining with the result of Theorem 3.1, we obtain Theorem 3.2. □

PROOF OF LEMMA 3.2. Obviously, the estimator $\hat{S}_n(t)$ is the same as $\hat{S}_n(t;\hat\beta_n,\hat\gamma_n)$. By repeating the proof of Theorem 3.2, we can obtain that if $|\beta - \beta^*| + |\gamma - \gamma^*| = o(a_n)$, then

$$\hat{S}_n(t;\beta,\gamma) - S(t)$$
$$= (\mathbf{P}_n - \mathbf{P})\left[ S_{T|Z^*}(t|Z^*) - S(t) \right.$$
$$- RI_{Y\leq t}S_{T|Z^*}(t|Z^*)e^{H_{T|Z^*}(Y|Z^*)+H_{C|Z^*}(Y|Z^*)}$$
$$\left. + S_{T|Z^*}(t|Z^*) \int_0^{t\wedge Y} e^{H_{T|Z^*}(u|Z^*)+H_{C|Z^*}(u|Z^*)} d_u H_{T|Z^*}(u|Z^*) \right]$$
$$- \mathbf{P}\{S_{T|Z^*}(t|Z^*)[B(\beta,\gamma,Z,t) - H_{T|Z^*}(t|Z^*)]\} + o_p\left(\frac{1}{\sqrt{n}}\right),$$

where we recall $B(\beta,\gamma,z,t) = \int_0^t \frac{d_u P(Y\leq u, R=1|(\beta'L,\gamma'L)=z)}{P(Y\geq u|(\beta'L,\gamma'L)=z)}$.

We especially choose $\gamma = \hat\gamma_n$ and $\beta = \hat\beta_n + \varepsilon_n v$ where $v$ is any constant vector on $R^{\dim(\beta^*)}$ with norm 1. After linearizing the $B(\beta,\gamma,Z,t)$ around $\beta = \beta^*$, $\gamma = \gamma^*$, we find that

$$\hat{S}_n(t;\hat\beta_n + \varepsilon_n v, \hat\gamma_n) - S(t)$$
$$= -\mathbf{P}\{S_{T|Z^*}(t|Z^*)[B(\beta^*,\gamma^*,Z^*,t) - H_{T|Z^*}(t|Z^*)]\}$$
$$- \varepsilon_n \mathbf{P}\{S_{T|Z^*}(t|Z^*)\nabla_\beta B(\beta^*,\gamma^*,Z^*,t)\}v + O_p\left(\frac{1}{\sqrt{n}}\right) + O(\varepsilon_n^2).$$

When one of the working models is correct,
$$-\mathbf{P}\{S_{T|Z^*}(t|Z^*)[B(\beta^*,\gamma^*,Z^*,t) - H_{T|Z^*}(t|Z^*)]\} = 0.$$

Moreover, $\hat{S}_n(t) - S(t) = O_p(\frac{1}{\sqrt{n}})$. Therefore,

$$\frac{\hat{S}_n(t;\hat\beta_n + \varepsilon_n v, \hat\gamma_n) - \hat{S}_n(t)}{\varepsilon_n} \xrightarrow{P} -\mathbf{P}\{S_{T|Z^*}(t|Z^*)\nabla_\beta B(\beta^*,\gamma^*,Z^*,t)\}v.$$

Similarly, for any constant vector $\tilde{v}$ in $R^{\dim(\gamma^*)}$ with norm 1,

$$\frac{\hat{S}_n(t;\hat{\beta}_n,\hat{\gamma}_n+\varepsilon_n\tilde{v})-\hat{S}_n(t)}{\varepsilon_n} \xrightarrow{P} -\mathbf{P}\{S_{T|Z^*}(t|Z^*)\nabla_\gamma B(\beta^*,\gamma^*,Z^*,t)\}\tilde{v}.$$

So the conclusions in the lemma hold. $\square$

**Acknowledgments.** This work is part of my Ph.D. dissertation advised by Professor Susan Murphy at the University of Michigan. I owe many thanks to her for numerous discussions and helpful comments. I also thank an Associate Editor and a referee for their valuable suggestions.

DEPARTMENT OF BIOSTATISTICS
UNIVERSITY OF NORTH CAROLINA
CHAPEL HILL, NORTH CAROLINA 27599-7420
USA
E-MAIL: dzeng@bios.unc.edu